\documentclass{amsproc}
\usepackage{amsfonts}

\theoremstyle{plain}

\newtheorem{example}{Example}

\newtheorem{lemma}{Lemma}

\newtheorem{proposition}{Proposition}

\numberwithin{equation}{section}

\input{tcilatex}

\begin{document}
\title[From an iteration formula to Poincar\'{e}'s Center Theorem]{From an
iteration formula to Poincar\'{e}'s Isochronous Center Theorem for
holomorphic vector fields}
\author{Guang Yuan Zhang}
\address{Department of Mathematical Sciences, Tsinghua University, Beijing
100084, The People's Republic of China}
\email{gyzhang@math.tsinghua.edu.cn\\
gyzhang@mail.tsinghua.edu.cn}
\date{}
\subjclass[2000]{32H50, 32M25, 37C27}
\keywords{ordinary differential equation, holomorphic differential equation}
\thanks{The author is supported by Chinese NSFC 10271063 and 10571009 }

\begin{abstract}
We first generalize a classical iteration formula for one variable
holomorphic mappings to a formula for higher dimensional holomorphic
mappings. Then, as an application, we give a short and intuitive proof of a
classical theorem, due to H. Poincar\'{e}, for the condition under which a
singularity of a holomorphic vector field is an isochronous center.
\end{abstract}

\maketitle

\section{An iteration formula}

Let $f$ be a holomorphic function germ at the origin in the complex plane
with $f(0)=0.$ For each positive integer $k,$ we denote by $f^{k}$ the $k$%
-th iteration of $f$ defined as $f^{1}=f,$ $f^{2}=f\circ f,$ ... , $%
f^{k}=f\circ f^{k-1}$ inductively, which is a well defined holomorphic
function germ at the origin.

Assume that $\lambda =f^{\prime }(0)$ is a primitive $m$-th root of unity,
say, $\lambda ^{m}=1$ but $\lambda ^{j}\neq 1$ for each positive integer $j$
with $j<m$. Then it is interesting that there exists a positive integer $r,$
such that the $m$-th iteration $f^{m}$ has a power series expansion 
\begin{equation*}
f^{m}(z)=z+a_{1}z^{rm+1}+a_{2}z^{rm+2}+\dots
\end{equation*}%
at the origin: all the terms of degrees from $2\ $to $rm$ vanish! This can
be proved by applying Rouch\'{e}'s theorem (see [\ref{Mi}]). In this section
we generalize this formula to germs of higher dimensional mappings by using
normal form theory. We denote by $\mathbb{C}^{n}$ the complex vector space
and by $\Delta ^{n}$ a ball in $\mathbb{C}^{n}$ centered at the origin.

\begin{proposition}[Iteration Formula]
\label{key}Let $f:\Delta ^{n}\rightarrow \mathbb{C}^{n}$ be a holomorphic
mapping given by 
\begin{equation*}
f(z)=\lambda z+o(|z|),\ z\in \Delta ^{n},
\end{equation*}%
where $\lambda $ is a primitive $m$-th root of unity. Then, in a
neighborhood of the origin, 
\begin{equation}
f^{m}(z)=z+o(|z|^{m}).  \label{ad0}
\end{equation}
\end{proposition}

In the proposition, $z=(z_{1},\dots ,z_{n})$ and the expression $o(|z|^{m})$
means that each component of the mapping $f^{m}(z)-z$ is a power series in $%
z_{1},\dots ,z_{n}$ consisting of terms of degree $>m.$

\begin{proof}
By the hypothesis and a fundamental result in the normal form theory (see [%
\ref{Ar}] or pages 84--85 in [\ref{AP}] for the proof), there exists a
biholomorphic transformation in the form of 
\begin{equation}
z=(z_{1},\dots ,z_{n})=h(x_{1,}\dots ,x_{n})=(x_{1},\dots ,x_{n})+o(|x|)
\label{ad-0}
\end{equation}%
of coordinate in a neighborhood of the origin such that each component $%
g_{j} $ of $g=h^{-1}\circ f\circ h=(g_{1},\dots ,g_{n})$ has a power series
expansion 
\begin{equation}
g_{j}(x_{1},\dots ,x_{n})=\lambda x_{j}+\sum c_{i_{1}\dots
i_{n}}^{j}x_{1}^{i_{1}}\dots x_{n}^{i_{n}}+o(|x|^{m}),\ \ j=1,\dots ,n,
\label{ad1}
\end{equation}%
in a neighborhood of the origin, where $x=(x_{1},\dots ,x_{n})$ and the sum
extends over all $n$-tuples $(i_{1},\dots ,i_{n})$ of nonnegative integers
with 
\begin{equation*}
2\leq i_{1}+\dots +i_{n}\leq m\mathrm{\ and\ }\lambda ^{i_{1}+\dots
+i_{n}}=\lambda .
\end{equation*}

On the other hand, since $\lambda $ is a primitive $m$-th root of unity, we
have 
\begin{equation*}
\lambda ^{i_{1}+\dots +i_{n}}\neq \lambda \ \mathrm{if\ }2\leq i_{1}+\dots
+i_{n}\leq m.
\end{equation*}%
Therefore, the sum in the equation (\ref{ad1}) vanishes, and then we have%
\begin{equation*}
g(x)=\lambda x+o(|x|^{m}),
\end{equation*}%
and then, considering that $\lambda ^{m}=1,$ we conclude that the $m$-th
iteration $g^{m}$ can be expressed as%
\begin{equation}
g^{m}(x)=x+o(|x|^{m}).  \label{f1}
\end{equation}

By (\ref{ad-0}) it is clear that $o(|h^{-1}(z)|^{m})=o(|z|^{m}),$ and then,
by (\ref{f1}) it is easy to see that 
\begin{eqnarray*}
&&f^{m}(z)=h\circ g^{m}\circ h^{-1}(z)=h(h^{-1}(z)+o(|h^{-1}(z)|^{m})) \\
&=&h(h^{-1}(z)+o(|z|^{m}))=h(h^{-1}(z))+o(|z|^{m}) \\
&=&z+o(|z|^{m}).
\end{eqnarray*}%
This completes the proof.
\end{proof}

\textbf{Acknowledgement }\textit{The author would like to give thanks to
Professor Meirong Zhang who told him that the normal form method is useful
for iteration problems five years ago and to Professor Xiaofeng Wang who
patiently taught him the fundamental result on normal forms which is used in
the above proof.}

\section{Poincar\'{e}'s condition for isochronous centers}

Consider an $n$-dimensional complex holomorphic system%
\begin{equation}
\dot{z}=F(z),\ \ z\in \Delta ^{n},  \label{1}
\end{equation}%
such that the origin $z=0$ is a singularity. This means that $F$ is a
holomorphic mapping from $\Delta ^{n}$ into $\mathbb{C}^{n}$ such that $%
F(0)=0$.

The origin is called a \emph{center }of the system if it has a punctured
neighborhood that is filled with periodic orbits, and called an \emph{%
isochronous} center if it has a punctured neighborhood that is filled with
periodic orbits of the same period. Here, and through out this paper, the 
\emph{period} of a periodic orbit means the smallest positive one.

The problem to find the condition so that a singularity of a system is a
center has a long history. The term \textit{center} was defined by H. Poincar%
\'{e}, while the research of center phenomena were started in 1673 when
Huygens studied the cycloidal pendulum (see [\ref{F}]). As an application of
the previous proposition, we shall present a short and intuitive proof of
the following classical theorem due to H. Poincar\'{e}.\medskip

\textbf{Isochronous Center Theorem. }\textit{If the Jacobian matrix }$%
F^{\prime }(0)$\textit{\ of }$F$\textit{\ at the origin\ equals to }$\omega
iI$\textit{\ for some real number }$\omega \neq 0,$\textit{\ where }$i=\sqrt{%
-1}$\textit{\ and }$I$\textit{\ is the unit matrix, then the origin\ is an
isochronous center with period }$2\pi /|\omega |.\medskip $

This result follows from Poincar\'{e}'s linearization theorem, which asserts
that the system (\ref{1}) is linearizable at the origin via a biholomorphic
transformation of coordinates, provided that the $n$-tuple $(\lambda
_{1},\dots ,\lambda _{n})$ of all eigenvalues of $F^{\prime }(0)$ is in the
Poincar\'{e} domain: the convex hull of these eigenvalues in the complex
plane does not contain the origin, and that there is no resonance: for any $%
n $-tuple $(i_{1},\dots ,i_{n})$ of nonnegative integers with $i_{1}+\dots
+i_{n}\geq 2$, $\lambda _{j}\neq \lambda _{1}^{i_{1}}\dots \lambda
_{n}^{i_{n}}$ for each $j=1,2,\dots ,n$ (see Chapter 5 in [\ref{Ar}] for the
details).

It is interesting that in the history of the study of central singularities,
special cases of the Isochronous Center Theorem have been rediscovered
several times. For example, when $n=1,$ it was rediscovered by Gregor [\ref%
{G1}] in 1958, Lukashevich [\ref{L}] in 1965, Brickman-Thomas [\ref{BT}] in
1977, Villarini [\ref{V}] in 1992, and Christopher-Devlin [\ref{CD}] in 1997
(see [\ref{JS}], [\ref{CRZ}], [\ref{Fr}], [\ref{H1}], [\ref{Pa}] and [\ref%
{Sa}] for other proofs and related topics for the case $n=1$). In 1998,
Needham-McAllister [\ref{NM}] rediscovered the result for two-dimensional
systems via the singularity theory of C. H. Briot and J. C. Bouquet. It
seems the approach in [\ref{NM}] applies to arbitrary dimensional case.

It is easy to see that a necessary condition so that the origin is an
isochronous center of the system (\ref{1}) is that all eigenvalues of $%
F^{\prime }(0)$ are pure imaginary with the same absolute value. The
converse fails in general.

\begin{example}
For the system 
\begin{equation*}
(\dot{x},\dot{y})=(ix,-iy+xy^{2}),\ (x,y)\in \mathbb{C}^{2},
\end{equation*}%
it is easy to verify that the corresponding flow is given by%
\begin{equation*}
\phi (t,(x,y))=(xe^{it},\frac{ye^{-it}}{1-txy}).
\end{equation*}%
Clearly, the origin is not a center of the system.
\end{example}

\section{Proof of the Isochronous Center Theorem}

The following result is well known (see [\ref{Zh3}] for a simple proof of a
more general version, where the singularity is just assumed to be isolated).

\begin{lemma}
\label{C1}If the Jacobian matrix $F^{\prime }(0)$ of $F$ at the origin is
invertible, then there exist a positive number $T_{0}$ and a ball $B$
centered at the origin, such that the system (\ref{1}) has no periodic orbit
in $B^{\ast }=B\backslash \{0\}$ with period less than $T_{0}.$
\end{lemma}

\begin{lemma}
\label{C2}There is a ball $B\subset \Delta ^{n}$ centered at the origin such
that the local flow $\phi (t,z)$ of (\ref{1}) is well defined and real
analytic on $[0,1]\times B,$ complex holomorphic with respect to $z,$ 
\begin{equation*}
\phi ([0,1]\times B)\subset \Delta ^{n},
\end{equation*}%
and, for each real number $\tau \in \lbrack 0,1],$ the Jacobian matrices $%
\Phi _{\tau }^{\prime }(0)$ and $F^{\prime }(0)$ of the time-$\tau $ map $%
\Phi _{\tau }$ and $F$ at the origin, respectively, satisfy%
\begin{equation*}
\Phi _{\tau }^{\prime }(0)=e^{\tau F^{\prime }(0)}.
\end{equation*}
\end{lemma}

The previous result is fundamental to the theory of holomorphic vector
fields. The time-$\tau $ map indicates the mapping $\Phi _{\tau
}:B\rightarrow \mathbb{C}^{n}$ given by 
\begin{equation*}
\Phi _{\tau }(z)=\phi (\tau ,z),z\in B,
\end{equation*}%
and the expression $e^{\tau F^{\prime }(0)}$ means the matrix $%
\sum_{k=0}^{\infty }\frac{\left( \tau F^{\prime }(0)\right) ^{n}}{n!}.$ If $%
F^{\prime }(0)=2\pi iI,$ for example, then $e^{\tau F^{\prime }(0)}=e^{2\pi
i\tau }I,$ where $I$ is the unit matrix.

\begin{proof}[\textbf{Proof of the Isochronous Center Theorem}]
Without loss of generality, assume $F^{\prime }(0)=2\pi iI.$ We shall show
that the origin is an isochronous center with period $1.$

Let $B$ be a ball centered at the origin that is determined by Lemma \ref{C2}%
. Let $\Phi _{1}:B\rightarrow \mathbb{C}^{n}$ be the time-$1$ map of the
local flow $\phi (t,z)$ of the system (\ref{1}). Then by Lemma \ref{C2} we
have $\Phi _{1}^{\prime }(0)=I,$ which implies that the equation%
\begin{equation}
\Phi _{1}(z)=z+o(|z|^{m}),\ z\in B,  \label{ad3}
\end{equation}%
holds for $m=1.$ But, we can show that this equation holds for all integers $%
m\geq 1$!

For any positive integer $m>1,$ consider the time-$\frac{\mathbf{1}}{m}$ map 
$\Phi _{\frac{1}{m}}$ of the local flow $\phi $. By Lemma \ref{C2}, we have $%
\Phi _{\frac{1}{m}}^{\prime }(0)=e^{\frac{2\pi i}{m}}I,$ and then by the
proposition, the $m$-th iteration $\Phi _{\frac{1}{m}}^{m}\ $of $\Phi _{%
\frac{1}{m}}$ can be expressed as%
\begin{equation*}
\Phi _{\frac{1}{m}}^{m}(z)=z+o(|z|^{m}),\ z\in B.
\end{equation*}%
Thus (\ref{ad3}) holds for all positive integers $m$. For 
\begin{equation*}
\Phi _{\frac{1}{m}}^{m}(z)=\phi (\frac{m}{m},z)=\phi (1,z)=\Phi _{1}(z),\
z\in B.
\end{equation*}

Since $\Phi _{1}$ is a holomorphic mapping on $B$, we have proved that $\Phi
_{1}=id_{B},$ the identity mapping on $B.$ Therefore, all orbits of the
system that intersect $B^{\ast }=B\backslash \{0\}$ are periodic, and then
the origin is a center.

Now, let us show that the origin is, in fact, an isochronous center with
period $1.$ We first show that for each $k\geq 2,$ the origin is an isolated
fixed point of the time-$\frac{\mathbf{1}}{k}$ map $\Phi _{\frac{1}{k}},$ of
the local flow. Otherwise, by the inverse function theorem, the Jacobian
determinant of the mapping $z\mapsto \Phi _{\frac{1}{k}}(z)-z$ vanishes at $%
z=0,$ and then $\Phi _{\frac{1}{k}}^{\prime }(0)$ must have an eigenvalue
equal to $1.$ But by Lemma \ref{C2}, $\Phi _{\frac{1}{k}}^{\prime }(0)=e^{%
\frac{2\pi i}{k}}I,$ and then we have $k=1$. Contradiction! Thus the origin
is an isolated fixed point of $\Phi _{\frac{1}{k}}.$ Therefore, for an
arbitrarily given integer $k_{0}\geq 2,$ there exists a neighborhood of the
origin in which the system has no periodic orbit of periods $\frac{1}{2}%
,\dots ,\frac{1}{k_{0}}.$

On the other hand, by the equation $\Phi _{1}=id_{B},$ it is clear that the
period of any periodic orbit of the system that intersects $B^{\ast }$ must
divide $1$, and then it equals to $\frac{\mathbf{1}}{k}$ for some positive
integer $k.$

Thus, there exists a punctured neighborhood of the origin in which each
periodic orbit of the system either has period $1$, or has period less than $%
\frac{1}{k_{0}}.$ Hence, by Lemma \ref{C1} and the arbitrariness of $k_{0}$,
we conclude that there exists a puncture neighborhood of the origin in which
all periodic orbits of the system has period $1.$ This completes the proof.
\end{proof}


\begin{thebibliography}{99}
\bibitem{} \label{Ar}Arnold, V.I. Geometrical methods in the theory of
ordinary differential equations, second edition, Translated by Joseph Sz\"{u}%
cs, English Translated Edited by Mark Levi, Springer-Verlag, 1988.

\bibitem{} \label{AP}Arrowsmith, D. K. \& Place, C. M., An introduction to
dynamical systems, Cambridge University Press, Cambridge, 1990.

\bibitem{} \label{BT} Brickman, L. \& Thomas, E. S., Conformal equivalence
of analytic flows, J. Differential Equations 25 (1977), no. 3, 310--324 (MR
56 \#5984 ).

\bibitem{} \label{JS} Chavarriga, J. \& Sabatini, M., A Survey of
isochronous centers, Qualitative theory of dynamical systems, 1 (1999),
1--70 (MR 2001c:34056).

\bibitem{} \label{CRZ} Cherkas, L. A., Romanovskii, V. G. \& \.{Z}o\l \c{a}%
dek, H., The centre conditions for a certain cubic system, Planar nonlinear
dynamical systems (Delft, 1995), Differential Equations Dynam. Systems 5
(1997), no. 3--4, 299--302 (MR\ 99i:34041).

\bibitem{} \label{CD} Christopher, C. J. \& Devlin, J., Isochronous centers
in planar polynomial systems, SIAM J. Math. Anal. 28 (1997), no. 1, 162--177
(MR 97k:34058).

\bibitem{} \label{F} Feigenbaum, L., The center of oscillation versus the
textbook writers of the early 18th century. From ancient omens to
statistical mechanics, 193--202, Acta Hist. Sci. Nat. Med., 39, Univ. Lib.
Copenhagen, Copenhagen, 1987 (MR 90m:01015).

\bibitem{} \label{Fr}Francoise, J.-P., Isochronous systems and perturbation
theory, Journal Nonlinear Math. Phys., Vol. 12, Supplement1 (2005), 315--326.

\bibitem{} \label{G1}Gregor, J., Dynamical systems with regular hand-side,
Pokroky Mat. Fys. Astronom. 3 (1958), 153--160 (Zbl 081.30802).

\bibitem{} \label{H1} Hajek, O., Notes on meromorphic dynamical systems,
I--III, Czechoslovak Math. J. 16 (1966), 14--40 (MR 33 \#2870a, 33 \#2870b,
33 \#2870c).

\bibitem{} \label{L} Lukashevich, N.A., Isochronicity of center for certain
systems of differential equations, Differ. Uravn. 1 (1965), 295--302 (MR 33
\#6023).

\bibitem{} \label{Mi} Milnor, J., Dynamics in One Complex Variable:
Introductory Lectures, Friedrick Vieweg \& Son, 2000

\bibitem{} \label{NM} Needham, D. J. \& McAllister, S., Centre families in
two-dimensional complex holomorphic dynamical systems, R. Soc. Lond. Proc.
Ser. A Math. Phys. Eng. Sci. 454 (1998), no. 1976, 2267--2278 (MR\
99d:34010).

\bibitem{} \label{Pa} Paluszny, M., On periodic solutions of polynomial ODEs
in the plane, \textit{J. Differential Equations} 53 (1984), no. 1, 24--29
(MR 86g:34054).

\bibitem{} \label{Sa} Sabatini, M., Dynamics of commuting systems on
two-dimensional manifolds, Ann. Mat. Pura Appl. (4) 173 (1997), 213--232 (MR
99f:34071).

\bibitem{} \label{V} Villarini, M., Regularity properties of the period
function near a center of a planar vector field, Nonlinear Anal. 19 (1992),
no. 8, 787--803 (MR 93j:34061).

\bibitem{} \label{Zh3}Zhang, G. Y., Fixed Point Indices and Invariant
Periodic Sets of Holomorphic Systems, to appear in Proc. Amer. Math. Soc.
\end{thebibliography}
\end{document}